\documentclass[12pt, reqno]{amsart}
\usepackage{amsmath}
\usepackage{amssymb}
\usepackage{epsfig}
\usepackage{graphicx}
\usepackage{color}
\usepackage{fullpage}
\usepackage{comment}
\definecolor{shadecolor}{gray}{0.875}
\usepackage{amscd}
\usepackage{stmaryrd}

\numberwithin{equation}{section}

\input xy
\xyoption{all}

\calclayout
\allowdisplaybreaks[3]

\theoremstyle{plain}
\newtheorem{prop}{Proposition}[section]

\theoremstyle{definition}
\newtheorem{defi}[prop]{Definition}

\newtheorem{conj}[prop]{Conjecture}

\newtheorem{rema}[prop]{Remark}

\newtheorem{exam}[prop]{Example}

\makeatother
\makeatletter

\author{Sho Tanimoto}
\address{Department of Mathematics, Faculty of Science, Kumamoto University, Kurokami 2-39-1 Kumamoto 860-8555 Japan}
\address{Priority Organization for Innovation and Excellence, Kumamoto University}
\email{stanimoto@kumamoto-u.ac.jp}

\title[Campana points and log Manin's conjecture]{Campana points, Height zeta functions, \\ and log Manin's conjecture}

\begin{document}
\date{\today}

\begin{abstract}
This is a report of the author's talk at RIMS workshop 2020 Problems and Prospects in Analytic Number Theory held online on Zoom.
We discuss a recent formulation of log Manin's conjecture for klt Campana points and an approach to this conjecture using the height zeta function method.
\end{abstract}

\maketitle

\section{Introduction}

One of fundamental tools in diophantine geometry is the notion of height functions and this height function measures the geometric and arithmetic complexities of rational points on an algebraic variety. These are crucial to various finiteness results in diophantine geometry such as Mordell-Weil theorem, Siegel's theorem, Mordell-Faltings' theorem, and so on. One of basic properties of height functions is the Northcott property which claims that for a height function associated to an ample divisor, the set of rational points whose height is less than $T$ is finite. Thus one may consider the counting function of rational points of bounded height, and one natural question is the asymptotic formula for such a counting function when $T$ goes to infinity.

Around the late 1980's, Yuri Manin and his collaborators proposed a general framework to understand this asymptotic formula in terms of geometric and arithmetic invariants of the underlying projective variety, and this leads to Manin's conjecture whose formulation is developed in a series of papers \cite{FMT89}, \cite{BM90}, \cite{Peyre}, \cite{BT98}, \cite{Peyre03}, \cite{Pey17}, and \cite{LST18}. One of fertile testing grounds for this conjecture is a class of equivariant compactifications of homogeneous spaces, and there are mainly two methods available, i.e., the method of mixing and the height zeta function method.

Mixing is a concept from ergodic theory, and this idea has been successfully used to prove equidistribution of rational points on homogeneous spaces acted by semi-simple groups (\cite{GMO08} and \cite{GO11}). The height zeta function method can be applied to a variety of equivariant compactifications of connected algebraic groups including, but not limited to, generalized flag varieties (\cite{FMT89}), toric varieties (\cite{BT96} and \cite{BT98}), equivariant compactifications of vector groups (\cite{CLT02}), wonderful compactifications of semi-simple groups of adjoint type (\cite{STBT07}), and biequivariant compactifications of unipotent groups (\cite{ST16}).

The height zeta function method also has its advantage to studying the counting problem of integral points associated to a reduced boundary divisor, and this has been implemented for equivariant compactifications of vector groups (\cite{CLT12}), toric varieties (\cite{CLT10toric}), wonderful compactifications of semi-simple groups of adjoint type (\cite{TBT13} and \cite{Chow19}), and biequivariant compactificaitons of the Heisenberg group (\cite{Xiao20}). These results suggest that there should be an analogous formulation of log Manin's conjecture for integral points, however certain subtleties of geometric and arithmetic nature prevent a general formulation of such a conjecture.

Campana and subsequently Abramovich proposed the notion of Campana points in \cite{Cam05} and \cite{Abr09} and this notion interpolates between rational points and integral points. The counting problem of Campana points has been originally featured in \cite{VVal11}, \cite{BVV12}, and \cite{VVal12}. Recently many mathematicians started to look at this problem and develop a series of results, attested by \cite{BY19}, \cite{PSTVA20}, \cite{PS20}, \cite{Xiao20}, and \cite{Str20}. In \cite{PSTVA20}, Pieropan, Smeets, V\'arilly-Alvarado, and the author initiated a systematic study of the counting problem for Campana points, and formulated a log Manin's conjecture for klt Campana points. Then we confirmed this conjecture for equivariant compactifications of vector groups using the height zeta function method for vector groups which is developed by Chambert-Loir and Tschinkel in \cite{CLT02} and \cite{CLT12}.

In this survey paper, we discuss the formulation of log Manin's conjecture for klt Campana points and applications of the height zeta function method to study this problem for various equivariant compactifications of connected algebraic groups.

\

Here is a plan of this paper: 
In Section~\ref{sec:height}, we review the notion of height functions. In Section~\ref{sec:campana}, we introduce two definitions of (weak) Campana points. In Section~\ref{sec:logmanin}, we discuss a formulation of log Manin's conjecture for klt Campana points. Finally in Section~\ref{sec:heightzeta}, we discuss the height zeta function method and its applications to equivariant compactifications of algebraic groups.

\

\noindent
{\bf Acknowledgements:}
The author would like to thank Marta Pieropan, Arne Smeets, and Tony V\'arilly-Alvarado for collaborations helping to shape his perspective on log Manin's conjecture for klt Campana points.
The author would like to thank the organizers of RIMS workshop Problems and Prospects in Analytic Number Theory for an opportunity to give a talk there. This work was supported by the Research Institute for Mathematical Sciences, an International Joint Usage/Research Center located in Kyoto University.

Sho Tanimoto was partially supported by Inamori Foundation, by JSPS KAKENHI Early-Career Scientists Grant number 19K14512, and by MEXT Japan, Leading Initiative for Excellent Young Researchers (LEADER).

\section{Height functions}
\label{sec:height}

In this section we review the notion of height functions and their basic properties. The main references are \cite{HS00} and \cite{CLT10}, and they include different treatments of height functions. In \cite{HS00} height functions are introduced using the machinery of Weil height machine and some basic properties of height functions such as the Northcott property are proved. In \cite{CLT10}, adelic metrizations are used to define height functions, and this definition is frequently used in the literature in Manin's conjecture. It is well-known that two definitions are essentially equivalent. See \cite{HS00} for more details. In this paper, we employ the definition of height functions using adelic metrizations described in \cite{CLT10}. 

\

Let us fix our notation: let $F$ be any number field and $\mathcal O_F$ be its ring of integers. We denote the set of places of $F$ by $\Omega_F$, the set of archimedean places by $\Omega_F^\infty$, and the set of non-archimedean places by $\Omega_F^{< \infty}$.
For any finite set $S \subset \Omega_F$ containing $\Omega_F^\infty$, $\mathcal O_{F, S}$ denotes the ring of $S$-integers. For each $v \in \Omega_F$, we denote the completion of $F$ with respect to $v$ by $F_v$. When $v$ is non-archimedean, we denote the ring of integers for $F_v$ by $\mathcal O_v$ with maximal ideal $\mathfrak m_v$ and residue field $k_v$ of size $q_v$. We denote the adele ring of $F$ by $\mathbb A_F$.

For each $v \in \Omega_F$, $F_v$ is a locally compact subgroup and it comes with a self-Haar measure $\mathrm dx_v = \mu_v$ which is normalized in a way Tate did in \cite{Tate}. We define the absolute value $|\cdot|_v$ on $F_v$ by requiring
\[
\mu_v(xB) = |x|_v \cdot \mu_v(B).
\]
This normalization satisfies the product formula, i.e., for any $x \in F^\times$, we have
\[
\prod_{v \in \Omega_F} |x|_v = 1
\]
See \cite{CLT10} for more details.

\

Let $F$ be a number field and $v \in \Omega_F$ be a place of $F$.
Let $U$ be an open set of $F_v^n$ in the analytic topology. A complex valued function on $U$ is smooth if it is $C^\infty$ when $v$ is archimedean and it is locally constant when $v$ is non-archimedean. This notion is local and extends to any $v$-adic analytic manifold.

Let $X$ be a smooth variety defined over $F_v$ and $L$ be a line bundle on $X$.
For each local point $x \in X(F_v)$, we denote the fiber of $L$ at $x$ by $L_x$.
\begin{defi}
A smooth metric on $L$ is a collection of metrics $\|\cdot\| : L_x(F_v) \rightarrow \mathbb R_{\geq 0}$ for all $x \in X(F_v)$ such that
\begin{itemize}
\item for $\ell \in L_x(F_v) \setminus \{0\}$, $\|\ell \| > 0$;
\item for any $a \in F_v$, $x \in X(F_v)$, and $\ell \in L_x(F_v)$, $\|a\ell \| = |a|_v \|\ell\|$, and;
\item for any open subset $U \subset X(F_v)$ and any non-vanishing section $\mathsf f \in \Gamma(U, L)$, the function $x \mapsto \|\mathsf f(x)\|$ is smooth.
\end{itemize}
\end{defi}

An integral model of a projective variety can be used to define a metric on it:
\begin{exam}
Let $X$ be a smooth projective variety defined over $F_v$ and $L$ be a line bundle on $X$ where $v$ is non-archimedean. Suppose that we have a flat projective $\mathcal O_v$-scheme $\mathcal X$ and a line bundle $\mathcal L$ on $\mathcal X$ extending $X$ and $L$. Let $x \in X(F_v) = \mathcal X(\mathcal O_v)$. Then we define a smooth metric on $L$ by insisting that for any $\ell \in L_x(F_v)$
\[
\|\ell\| \leq 1 \iff \ell \in \mathcal L_x(\mathcal O_v).
\]
This metric is called as the induced metric by an integral model $(\mathcal X, \mathcal L)$.
\end{exam}

Next we define adelic metrizations on a smooth projective variety defined over a number field $F$.
\begin{defi}
Let $X$ be a smooth projective variety defined over a number field $F$ and $L$ be a line bundle on $X$. An adelic metrization on $L$ is a collection of $v$-adic smooth metrics $\{\|\cdot\|_v\}_{v \in \Omega_F}$ on $X$ such that there exist a finite set $S$ of places including $\Omega_F^\infty$, a flat $\mathcal O_{F, S}$-projective model $\mathcal X$, and a line bundle $\mathcal L$ on $\mathcal X$ extending $X$ and $L$ such that for any $v \not\in S$, the metric $\|\cdot\|_v$ is induced by $(\mathcal X, \mathcal L)$. Note that two integral models are isomorphic outside of finitely many places so that two adelic metrizations differ only at finitely many places.
\end{defi}

Finally we define the notion of height functions:
\begin{defi}
Let $X$ be a smooth projective variety defined over a number field $F$ and $\mathcal L = (L, \{\|\cdot \|_v\})$ be an adelically metrized line bundle on $X$.
For each rational point $x \in X(F)$, choose $\ell \in L_x(F)$ and we define the height function $\mathsf H_{\mathcal L} : X(F) \to \mathbb R_{\geq 0}$ by
\[
\mathsf H_{\mathcal L}(x) = \prod_{v \in \Omega_F} \|\ell\|_v^{-1}.
\]
This is well-defined due to the product formula mentioned above.
\end{defi}

Here is an example of height functions for the projective space:
\begin{exam}
Let $X = \mathbb P^n$ and $L = \mathcal O_X(1)$. We consider the standard integral model $\mathcal X = \mathbb P^n_{\mathcal O_F}$. For each non-archimedean place $v \in \Omega_F$, we let $\|\cdot\|_v$ be the metric at $v$ induced by $\mathcal X$. For any archimedean place $v$, we define a smooth metric at $v$ by insisting
\[
\|\ell(x) \|_v = \frac{|\ell(x)|_v}{\sqrt{\sum_{i = 0}^n |x_i|_v^2}},
\]
where $x = (x_0: \cdots : x_n) \in X(F)$ and $\ell \in H^0(X, \mathcal O_X(1))$. Then it is an easy exercise to prove that the height function associated to $L$ with this adelic metrization is given by
\[
\mathsf H(x) = \prod_{v \in \Omega_F^{< \infty}} \max \{ |x_0|_v, \cdots, |x_n|_v\} \prod_{v \in \Omega_F^\infty} \sqrt{|x_0|_v^2 + \cdots + |x_n|_v^2}.
\]
When $F = \mathbb Q$, we may assume that $x_i$'s are integers and $\gcd (x_0, \cdots, x_n) = 1$. In this situation, the above formula reduces to
\[
\mathsf H(x) = \sqrt{|x_0|_\infty^2 + \cdots + |x_n|_\infty^2}.
\]
\end{exam}

Let us mention a few basic properties of height functions:
\begin{prop}
Let $X$ be a smooth projective variety defined over a number field $F$ and $\mathcal L = (L, \{\|\cdot\|_v\})$ be an adelically metrized line bundle on $X$.
Then the following statements are true:
\begin{itemize}
\item Let $\mathcal L'$ be another adelically metrized line bundle associated to $L$. Then there exist positive constants $C_1 \leq C_2$ such that for any $x \in X(F)$, we have
\[
C_1 \mathsf H_{\mathcal L'}(x) \leq \mathsf H_{\mathcal L}(x) \leq C_2 \mathsf H_{\mathcal L'}(x);
\]
\item Let $B$ be the base locus of the complete linear series $|L|$. Then there exists a positive constant $C > 0$ such that for any $x \in (X \setminus B)(F)$, we have
\[
\mathsf H_{\mathcal L}(x) \geq C;
\]
\item When $L$ is ample, for any real number $T > 0$ the set
\[
\{x \in X(F) \, \mid \, \mathsf H_{\mathcal L}(x)\ \leq T \}
\]
is a finite set.
\end{itemize}
\end{prop}

The last property is called as the Northcott property which is fundamental in diophantine geometry and it is also foudational for Manin's conjecture. For more details, see \cite{HS00}.

\section{Campana points}
\label{sec:campana}

In this section, we review the notion of Campana points. Campana points were originally considered by Campana for curves in \cite{Cam05}, and its higher dimensional analogue was explored by Abramovich in \cite{Abr09}. One may consider Campana points as integral points on Campana orbifolds developed by again Campana himself:

\begin{defi}
Let $F$ be an arbitrary field and $X$ be a smooth projective variety defined over $F$. Let $D_\epsilon = \sum_{\alpha \in \mathcal A} \epsilon_\alpha D_\alpha$ be an effective $\mathbb Q$-divisor on $X$ with $D_\alpha$'s irreducible and distinct. We say $(X, D_\epsilon)$ is a Campana orbifold if the following statements are true:
\begin{itemize}
\item For any $\alpha \in \mathcal A$, a non-negative rational number $\epsilon_\alpha$ takes the form of
\[
1- \frac{1}{m_\alpha},
\]
where $m_\alpha$ is a positive integer or $+\infty$;
\item the reduced divisor $D = \sum_{\alpha \in \mathcal A} D_\alpha$ is a strict normal crossings divisor.
\end{itemize}
We say a Campana orbifold $(X, D_\epsilon)$ is Fano if $-(K_X + D_\epsilon)$ is ample.
\end{defi}

Let $(X, D_\epsilon)$ be a Campana orbifold. Then $(X, D_\epsilon)$ is a divisorial log terminal (dlt for short) pair in the sense of birational geometry. When $\epsilon_\alpha < 1$ for any $\alpha$, $(X, D_\epsilon)$ is a kawamata log terminal (klt for short) pair. See \cite{KM98} for the definitions and their basic properties. We say a Campana orbifold $(X, D_\epsilon)$ is klt if $\epsilon_\alpha < 1$ for every $\alpha \in \mathcal A$.

\

To define the notion of Campana points, one needs to fix an integral model of a Campana orbifold.
Let $(X, D_\epsilon)$ be a Campana orbifold defined over a number field $F$ with $D_\epsilon = \sum_{\alpha \in \mathcal A} \epsilon_\alpha D_\alpha$. Let $S$ be a finite set of places including all archimedean places. A good integral model of $(X, D_\epsilon)$ away from $S$ is a flat projective $\mathcal O_{F, S}$-scheme $\mathcal X$ such that $\mathcal X$ is extending $X$ and $\mathcal X$ is regular. Let $\mathcal D_\alpha$ be the Zariski closure of $D_\alpha$ in $\mathcal X$ and let $\mathcal D_\epsilon = \sum_{\alpha \in \mathcal A} \epsilon_\alpha \mathcal D_\alpha$.

Let us fix a good integral model of a Campana orbifold $(X, D_\epsilon)$ as above. Let $\mathcal A_\epsilon = \{\alpha \in \mathcal A \, | \, \epsilon_\alpha \neq 0\}$. We set $X^\circ = X\setminus \cup_{\alpha \in \mathcal A_\epsilon} D_\alpha$. Let $P \in X^\circ(F)$ be a rational point and $v \not\in S$ be a non-archimedean place of $F$. Then we may consider $P$ as an $\mathcal O_v$-point $\mathcal P_v \in \mathcal X(\mathcal O_v)$ by valuative criterion for properness. Since $\mathcal P_v \not\subset \mathcal D_\alpha$ for any $\alpha \in \mathcal A_\epsilon$, the pullback of $\mathcal D_\alpha$ via $\mathcal P_v$ defines an ideal in $\mathcal O_v$. We denote its colength by $n_v(\mathcal D_\alpha, P)$. When $P \in D_\alpha$ for some $\alpha \in \mathcal A_\epsilon$, we formally set $n_v(\mathcal D_\alpha, P) = + \infty$. The total intersection number is given by
\[
n_v(\mathcal D_\epsilon, P) = \sum_{\alpha \in \mathcal A_\epsilon} \epsilon_\alpha n_v(\mathcal D_\alpha, P).
\]

Now we are ready to define two notions of Campana points:
\begin{defi}
We say $P \in X(F)$ is a weak Campana $\mathcal O_{F, S}$-point on $(\mathcal X, \mathcal D_\epsilon)$ if the following statements are true:
\begin{itemize}
\item we have $P \in (\mathcal X \setminus \cup_{\epsilon_\alpha = 1}\mathcal D_\alpha)(\mathcal O_{F, S})$, and;
\item for $v \not\in S$, if $n_v(\mathcal D_\epsilon, P) > 0$, then
\[
n_v(\mathcal D_\epsilon, P) \leq \left (\sum_{\alpha \in \mathcal A_\epsilon} n_v(\mathcal D_\alpha, P) \right)-1.
\]
\end{itemize}
\end{defi}
We denote the set of weak Campana $\mathcal O_{F, S}$-points by $(\mathcal X, \mathcal D_\epsilon)_w(\mathcal O_{F, S})$.

\begin{defi}
We say $P \in X(F)$ is a Campana $\mathcal O_{F, S}$-point on $(\mathcal X, \mathcal D_\epsilon)$ if the following statements are true:
\begin{itemize}
\item we have $P \in (\mathcal X \setminus \cup_{\epsilon_\alpha = 1}\mathcal D_\alpha)(\mathcal O_{F, S})$, and;
\item for $v \not\in S$ and for all $\alpha \in \mathcal A_\epsilon$ with $\epsilon_\alpha < 1$ and $n_v(\mathcal D_\alpha, P) > 0$, we have
\[
n_v(\mathcal D_\alpha, P) \geq m_\alpha,
\]
where $\epsilon_\alpha = 1-\frac{1}{m_\alpha}$.
\end{itemize}
A Campana $\mathcal O_{F, S}$-point is klt when the underlying Campana orbifold is a klt pair.
\end{defi}

We denote the set of Campana $\mathcal O_{F, S}$-points by $(\mathcal X, \mathcal D_\epsilon)(\mathcal O_{F, S})$. Then we have the following inclusions:
\[
\mathcal X^\circ(\mathcal O_{F, S}) \subset (\mathcal X, \mathcal D_\epsilon)(\mathcal O_{F, S}) \subset (\mathcal X, \mathcal D_\epsilon)_w(\mathcal O_{F, S}) \subset X(F),
\]
where $\mathcal X^\circ = \mathcal X \setminus (\cup_{\alpha \in \mathcal A_\epsilon}\mathcal D_\alpha)$. When $\epsilon_\alpha = 0$ for all $\alpha \in \mathcal A$, the rightmost two inclusions are equalities. When $\epsilon_\alpha = 1$ for all $\alpha \in \mathcal A_\epsilon$, the leftmost two inclusions are equalities.

Here is an example of klt Campana points:
\begin{exam}
For simplicity, let us assume that $F = \mathbb Q$ and $S = \{\infty\}$. Let $X = \mathbb P^n$ and $H = V(x_0)$ be a hyperplane. Let $m$ be a positive integer and $\epsilon = 1-1/m$. We define 
\[
D_\epsilon = \epsilon H.
\]
We consider the standard integral model of $X$. Then a rational point $x = (x_0 : \cdots: x_n) \in X(\mathbb Q)$ with $x_i \in \mathbb Z$ and $\gcd(x_0, \cdots, x_n) = 1$ is a Campana $\mathbb Z$-point if $x_0 = 0$ or $x_0 \neq 0$ and the following statement is true: for any prime number $p$ we have
\[
p \mid x_0 \implies p^m \mid x_0.
\]
Any non-zero integer with this property is said to be $m$-full. When $m = 2$, it is said to be squarefull.
\end{exam}

\section{Log Manin's conjecture}
\label{sec:logmanin}

Let $X$ be a smooth projective variety defined over a number field $F$ and $\mathcal L = (L, \{\|\cdot\|_v\})$ be an adelically metrized line bundle on $X$. We consider the associated height function
\[
\mathsf H_{\mathcal L} : X(F) \to \mathbb R_{>0}.
\]
When $L$ is ample, this height function satisfies the Northcott property so that for any subset $Q \subset X(F)$ and any positive real number $T>0$ one may define the counting function
\[
N(Q, \mathcal L, T) = \# \{P \in Q \, | \, \mathsf H_{\mathcal L}(P) \leq T\}.
\]
Manin's conjecture predicts the asymptotic formula of the above function for an appropriate $Q$, and a natural question is to extend this conjecture to integral points and Campana points. In \cite{PSTVA20}, Pieropan, Smeets, V\'arilly-Alvarado and the author formulated this log version of Manin's conjecture when the underlying Campana orbifold is a klt log Fano pair. In this section, we review a general formulation of this log Manin's conjecture.

\subsection{Two birational invariants}

Let $X$ be a smooth projective variety defined over a field $F$. Let $D_1, D_2$ are $\mathbb Q$-divisors on $X$. We say $D_1$ and $D_2$ are numerically equivalent if for any curve $C \subset X$, we have $D_1.C = D_2.C$. In this case we write $D_1 \equiv D_2$. We define the space of $\mathbb Q$-divisors up to numerical equivalence as
\[
N^1(X)_{\mathbb Q} = \{D : \text{$\mathbb Q$-divisors}\}/\equiv.
\]
We set $N^1(X) := N^1(X)_{\mathbb Q} \otimes_{\mathbb Q}\mathbb R$.
Then we define the cone of pseudo-effective divisors by
\[
\overline{\mathrm{Eff}}^1(X) := \overline{\text{the cone of effective $\mathbb Q$-divisors}} \subset N^1(X).
\]

Now we are ready to introduce two birational invariants which play central roles in Manin's conjecture:
\begin{defi}
Let $(X, D_\epsilon)$ be a klt Campana orbifold defined over a field $F$ and $L$ be an ample $\mathbb Q$-divisor on $X$. We define the Fujita invariant or $a$-invariant by
\[
a(X, D_\epsilon, L) := \inf \{t \in \mathbb R \, |\, tL + K_X + D_\epsilon \in \overline{\mathrm{Eff}}^1(X)\}.
\]

Next assume that $a(X, D_\epsilon, L) >0$. Then we define the $b$-invariant by
\begin{align*}
b(F, X, D_\epsilon, L) := & \text{ codimension of the minimal face of $\overline{\mathrm{Eff}}^1(X)$}\\ & \text{ containing $a(X, D_\epsilon, L)L + K_X + D_\epsilon$.}
\end{align*}
It is explained in \cite[Section 3.6.2]{PSTVA20} that these invariants are birational invariants.
\end{defi}

\begin{exam}
Let $(X, D_\epsilon)$ be a klt Fano orbifold defined over a field $F$ and $L = -(K_X + D_\epsilon)$.
Then we have
\[
a(X, D_\epsilon, L) = 1, \quad b(F, X, D_\epsilon, L) = \rho (X) = \dim N^1(X).
\]
\end{exam}

\subsection{Thin exceptional sets}

The notion of thin sets has been explored by Serre to study Galois inverse problem, and it is also fundamental to Manin's conjecture. Let us give the definition of thin sets for Campana points:
\begin{defi}
Let $(X, D_\epsilon)$ be a klt Campana orbifold defined over a number field $F$. Let $S$ be a finite set of places of $F$ including $\Omega_F^\infty$ and we fix a good integral model away from $S$ $\mathcal X \to \mathrm{Spec} \, \mathcal O_{F, S}$.

A type I thin set is a set of the form
\[
V(F) \cap (\mathcal X, \mathcal D_\epsilon)(\mathcal O_{F, S}),
\]
where $V \subset X$ is a proper closed subset of $X$.

A type II thin set is a set of the form
\[
f(Y(F))\cap (\mathcal X, \mathcal D_\epsilon)(\mathcal O_{F, S}), 
\]
where $f : Y \to X$ is a dominant generically finite morphism of degree $\geq 2$ defined over $F$ with $Y$ integral.

A thin set is any subset of a finite union of type I and type II thin sets.
\end{defi}

Here is an example of thin sets:
\begin{exam}
Let $X = \mathbb P^1$ with $D_\epsilon = 0$ defined over a number field $F$.
We consider the morphism
\[
f : \mathbb P^1 \to \mathbb P^1, \, (x_0:x_1) \mapsto (x_0^d:x_1^d)
\]
with $d \geq 2$.
Then $f(X(F)) \subset X(F)$ is a thin set.

\end{exam}

\subsection{Log Manin's conjecture for klt Campana points}

Finally we state log Manin's conjecture for klt Campana points:

\begin{conj}[Log Manin's conjecture for klt Campana points]
\label{conj:main}
Let $(X, D_\epsilon)$ be a klt Fano orbifold defined over a number field $F$ and $\mathcal L = (L, \{\|\cdot\|_v\})$ be an adelically metrized ample line bundle. Assume that $(\mathcal X, \mathcal D_\epsilon)(\mathcal O_{F, S})$ is not thin. Then there exists a thin set $Z \subset (\mathcal X, \mathcal D_\epsilon)(\mathcal O_{F, S})$ such that 
\[
N((\mathcal X, \mathcal D_\epsilon)(\mathcal O_{F, S}) \setminus Z, \mathcal L, T) \sim c(F, \mathcal X, \mathcal D_\epsilon, \mathcal L, Z)T^{a(X, D_\epsilon, L)} (\log T)^{b(F, X, D_\epsilon, L)-1},
\]
as $T \to \infty$. Here the leading constant $c(F, \mathcal X, \mathcal D_\epsilon, \mathcal L, Z)$ is analogous to Peyre's constant developed in \cite{Peyre} and \cite{BT98} and its definition is given in \cite[Section 3.3]{PSTVA20}.
\end{conj}

\begin{rema}
For a smooth geometrically rationally connected projective variety $X$ defined over a number field $F$, it is expected that $X(F)$ is not thin as soon as there is a rational point. Indeed, Colliot-Th\'etl\`ene's conjecture predicts that the set of rational points is dense in the Brauer-Manin set, and this implies that $X$ satisfies weak weak approximation. It is known that weak weak approximation property implies non-thinness of the set of rational points. The corresponding statement for klt Campana points, i.e., weak weak approximation for klt Campana sets implies non-thiness of the set of klt Campana points is established in \cite{NS20}. So it is natural to expect that the assumption of Conjecture~\ref{conj:main} is true as long as there is a klt Campana point.
\end{rema}

\begin{rema}
It is well-documented in the case of rational points that in Conjecture~\ref{conj:main} it is important to remove the contribution of a thin set $Z$ from the counting function. There is a series of papers (\cite{LTDuke}, \cite{Sen17}, and \cite{LST18}) studying birational geometry of thin exceptional subsets for rational points. In \cite{LST18}, Lehmann, Sengupta, and the author proposed a conjectural description of thin exceptional subsets and proved that it is indeed a thin set using the minimal model program and the boundedness of singular Fano varieties. It would be interesting to perform a similar study for klt Campana points.
\end{rema}

Conjecture~\ref{conj:main} is known in the following cases:
\begin{itemize}
\item projective space with a boundary being the union of hyperplanes (\cite{VVal11}, \cite{VVal12}, \cite{BVV12}, and \cite{BY19});
\item equivariant compactifications of vector groups (\cite{PSTVA20});
\item toric varieties defined over $\mathbb Q$ (\cite{PS20}) and;
\item biequivaraint compactifications of the Heisenberg group (\cite{Xiao20}).
\end{itemize}

One can also consider a similar counting problem for weak Campana points, however this problem is much harder than Conjecture~\ref{conj:main}. At the moment of writing this paper, we do not know how one should formulate a log Manin's conjecture for weak Campana points, but \cite{Str20} takes the first step towards to this problem.

\section{Height zeta functions}
\label{sec:heightzeta}

Let $F$ be a number field and $G$ be a connected linear algebraic group defined over $F$.
Let $X$ be a smooth projective equivariant compacitification of $G$, i.e., $X$ contains $G$ as a Zariski open subset, and the right action of $G$ extends to $X$. In this situation, the boundary 
\[
D = X \setminus G = \bigcup_{\alpha \in \mathcal A} D_\alpha
\]
is a divisor where each $D_\alpha$ is an irreducible component. After applying an equivariant resolution, we may assume that $D = \sum_{\alpha \in \mathcal A}D_\alpha$ is a divisor with strict normal crossings. We also fix an adelic metriazation for $\mathcal O(D_\alpha)$ for every $\alpha \in \mathcal A$.

For each $\alpha \in \mathcal A$, we choose $m_\alpha$ which is a positive integer or $+ \infty$, and set $\epsilon_\alpha = 1-\frac{1}{m_\alpha}$. We consider
\[
D_\epsilon = \sum_{\alpha \in \mathcal A}\epsilon_\alpha D_\alpha,
\]
and $(X, D_\epsilon)$ is a Campana orbifold.
Let us fix a finite set $S$ of places including $\Omega_F^\infty$ and a good integral model $\mathcal X$ away from $S$ extending $X$.
When $-(K_X + D_\epsilon)$ is ample (or more generally big), it is natural to consider the counting problem of $\mathcal O_{F, S}$-Campana points on $(\mathcal X, \mathcal D_\epsilon)$. There is a general approach to this problem which is called as the height zeta function method.

Let $\mathrm{Pic}(X)^G$ be the Picard group of $G$-linearlized line bundles on $X$ up to isomorphisms. (If the reader is not familiar with $G$-linearlizations, she/he may ignore this term for now.) After tensoring by $\mathbb Q$, boundary components $D_\alpha$ form a basis for $\mathrm{Pic}(X)^G_{\mathbb Q}$. 
We choose a section
\[
\mathsf f_\alpha \in H^0(X, \mathcal O(D_\alpha)),
\]
corresponding to $D_\alpha$. Then we define a local height pairing: for any place $v\in \Omega_F$,
\[
\mathsf H_v : G(F_v) \times \mathrm{Pic}(X)^G_{\mathbb C} \to \mathbb C^\times, \left(g_v, \sum_{\alpha \in \mathcal A} s_\alpha D_\alpha\right) \mapsto \prod_{\alpha \in \mathcal A} \|\mathsf f_\alpha (g_v)\|_v^{-s_\alpha}.
\]
Using this local height pairing, we define the global height pairing as the Euler product:
\[
\mathsf H := \prod_{v \in \Omega_F} \mathsf H_v : G(\mathbb A_F) \times \mathrm{Pic}(X)^G_{\mathbb C} \to \mathbb C^\times.
\]

Applying the definition of Campana points to local points, for each $v \not\in S$, one can define the Campana set
\[
(\mathcal X, \mathcal D_\epsilon)(\mathcal O_v) \subset X(F_v).
\]
We set 
\[
G(F_v)_\epsilon = G(F_v) \cap (\mathcal X, \mathcal D_\epsilon)(\mathcal O_v),
\]
and let $\delta_{\epsilon, v}(g_v)$ be the characteristic function of $G(F_v)_\epsilon$ on $G(F_v)$.
When $v \in S$, we set $\delta_{\epsilon, v}\equiv 1$ and define $\delta_\epsilon$ as the Euler product:
\[
\delta_\epsilon = \prod_{v \in \Omega_F} \delta_{\epsilon, v} : G(\mathbb A_F) \to \mathbb R_{\geq 0}.
\]

For $g \in G(\mathbb A_F)$ and $\bold s \in \mathrm{Pic}(X)^G_{\mathbb C}$, we define the height zeta function by
\[
\mathsf Z(g, \bold s) := \sum_{\gamma \in G(F)} \mathsf H(\gamma g, \bold s)^{-1}\delta_\epsilon(\gamma g)
\]
When $\Re(\bold s)$ is sufficiently large, this function converges to a continuous function in $g \in G(F)\backslash G(\mathbb A_F)$ and a holomorphic function in $\bold s \in \mathrm{Pic}(X)^G_{\mathbb C}$.

A relation of this height zeta function to log Manin's conjecture is given by Tauberian theorem. Indeed, if one can prove that for an ample (or big) line bundle $L$, $\mathsf Z(\mathrm{id}, sL)$ admits a meromorphic continuation to a half plane $\Re(s) \geq a$ with a unique pole at $s = a$ of order $b$ with $a >0$ a positive real number, then one can conclude
\[
N(G(F)_\epsilon, \mathcal L, T) \sim c T^{a} (\log T)^{b-1},
\]
where $c$ is a positive constant related to the leading constant of $\mathsf Z(\mathrm{id}, sL)$ at $s = a$. Thus our goal is reduced to obtain a meromorphic continuation of $\mathsf Z(\mathrm{id}, \bold s)$.

To this end, for $\bold s \gg 0$, one can prove that
\[
\mathsf Z(g, \bold s) \in \mathsf L^2(G(F)\backslash G(\mathbb A_F)),
\]
thus one may apply spectral decomposition of this Hilbert space to $\mathsf Z(g, \bold s)$ and use this spectral decomposition to obtain a meromorphic continuation.

This program has been pioneered mainly by Tschinkel and his collaborators, and has been carried out in the following cases:
\begin{itemize}
\item rational points on toric varieties (\cite{BT96}, \cite{BTtoric});
\item rational points on equivariant compacitifications of vector groups (\cite{CLT02});
\item rational points on wonderful compactifications of semi-simple groups of adjoint type (\cite{STBT07});
\item rational points on biequivariant compactifications of unipotent groups (\cite{ST16});
\item integral points on equivariant compacitificaitons of vector groups (\cite{CLT12});
\item integral points on toric varieties (\cite{CLT10toric});
\item integral points on wonderful compactificaitons of semi-simple groups of adjoint type (\cite{TBT13} and \cite{Chow19});
\item Campana points on equivariant compacitifications of vector groups (\cite{PSTVA20});
\item Campana points on biequivariant compacitifications of the Heisenberg group (\cite{Xiao20}), and;
\item weak Campana points on certain toric varieties (\cite{Str20}).
\end{itemize}
It would be interesting to explore Campana points on other algebraic groups. In particular, the treatment of integral points on toric varieties (\cite{CLT10toric}) is known to be incomplete, and there is some technical issue on this paper. It would be interesting to apply the height zeta function method to klt Campana points on toric varieties and see whether we have a similar issue.

Finally for the readers who are interested in working examples of this program, we recommend them to consult \cite[Interlude I]{PSTVA20}.

\nocite{*}
\bibliographystyle{alpha}
\bibliography{RIMS}

\def\cprime{$'$}
\begin{thebibliography}{PSTVA20}

\bibitem[Abr09]{Abr09}
Dan Abramovich.
\newblock Birational geometry for number theorists.
\newblock In {\em Arithmetic geometry}, volume~8 of {\em Clay Math. Proc.},
  pages 335--373. Amer. Math. Soc., Providence, RI, 2009.

\bibitem[BM90]{BM90}
V.~V. Batyrev and Yu.~I. Manin.
\newblock Sur le nombre des points rationnels de hauteur born\'{e} des
  vari\'{e}t\'{e}s alg\'{e}briques.
\newblock {\em Math. Ann.}, 286(1-3):27--43, 1990.

\bibitem[BT96]{BT96}
V.~Batyrev and Yu. Tschinkel.
\newblock Height zeta functions of toric varieties.
\newblock volume~82, pages 3220--3239. 1996.
\newblock Algebraic geometry, 5.

\bibitem[BT98a]{BT98}
V.~V. Batyrev and Y.~Tschinkel.
\newblock Tamagawa numbers of polarized algebraic varieties.
\newblock Number 251, pages 299--340. 1998.
\newblock Nombre et r\'{e}partition de points de hauteur born\'{e}e (Paris,
  1996).

\bibitem[BT98b]{BTtoric}
Victor~V. Batyrev and Yuri Tschinkel.
\newblock Manin's conjecture for toric varieties.
\newblock {\em J. Algebraic Geom.}, 7(1):15--53, 1998.

\bibitem[BVV12]{BVV12}
T.~D. Browning and K.~Van~Valckenborgh.
\newblock Sums of three squareful numbers.
\newblock {\em Exp. Math.}, 21(2):204--211, 2012.

\bibitem[BY20]{BY19}
T.~D. Browning and S.~Yamagishi.
\newblock Arithmetic of higher-dimensional orbifolds and a mixed {W}aring
  problem.
\newblock {\em Math. Z.}, 2020.
\newblock to appear.

\bibitem[Cam05]{Cam05}
Fr\'{e}d\'{e}ric Campana.
\newblock Fibres multiples sur les surfaces: aspects geom\'{e}triques,
  hyperboliques et arithm\'{e}tiques.
\newblock {\em Manuscripta Math.}, 117(4):429--461, 2005.

\bibitem[Cho19]{Chow19}
D.~Chow.
\newblock The {D}istribution of {I}ntegral {P}oints on the{W}onderful
  {C}ompactifications by {H}eight.
\newblock submitted, 2019.

\bibitem[CLT02]{CLT02}
Antoine Chambert-Loir and Yuri Tschinkel.
\newblock On the distribution of points of bounded height on equivariant
  compactifications of vector groups.
\newblock {\em Invent. Math.}, 148(2):421--452, 2002.

\bibitem[CLT10a]{CLT10}
Antoine Chambert-Loir and Yuri Tschinkel.
\newblock Igusa integrals and volume asymptotics in analytic and adelic
  geometry.
\newblock {\em Confluentes Math.}, 2(3):351--429, 2010.

\bibitem[CLT10b]{CLT10toric}
Antoine Chambert-Loir and Yuri Tschinkel.
\newblock Integral points of bounded height on toric varieties.
\newblock preprint, 2010.

\bibitem[CLT12]{CLT12}
Antoine Chambert-Loir and Yuri Tschinkel.
\newblock Integral points of bounded height on partial equivariant
  compactifications of vector groups.
\newblock {\em Duke Math. J.}, 161(15):2799--2836, 2012.

\bibitem[FMT89]{FMT89}
Jens Franke, Yuri~I. Manin, and Yuri Tschinkel.
\newblock Rational points of bounded height on {F}ano varieties.
\newblock {\em Invent. Math.}, 95(2):421--435, 1989.

\bibitem[GMO08]{GMO08}
Alex Gorodnik, Fran\c{c}ois Maucourant, and Hee Oh.
\newblock Manin's and {P}eyre's conjectures on rational points and adelic
  mixing.
\newblock {\em Ann. Sci. \'{E}c. Norm. Sup\'{e}r. (4)}, 41(3):383--435, 2008.

\bibitem[GO11]{GO11}
Alex Gorodnik and Hee Oh.
\newblock Rational points on homogeneous varieties and equidistribution of
  adelic periods.
\newblock {\em Geom. Funct. Anal.}, 21(2):319--392, 2011.
\newblock With an appendix by Mikhail Borovoi.

\bibitem[HS00]{HS00}
Marc Hindry and Joseph~H. Silverman.
\newblock {\em Diophantine geometry}, volume 201 of {\em Graduate Texts in
  Mathematics}.
\newblock Springer-Verlag, New York, 2000.
\newblock An introduction.

\bibitem[KM98]{KM98}
J.~Koll\'{a}r and Sh. Mori.
\newblock {\em Birational geometry of algebraic varieties}, volume 134 of {\em
  Cambridge Tracts in Mathematics}.
\newblock Cambridge University Press, Cambridge, 1998.
\newblock With the collaboration of C. H. Clemens and A. Corti, Translated from
  the 1998 Japanese original.

\bibitem[LST18]{LST18}
B.~Lehmann, A.~K. Sengupta, and S.~Tanimoto.
\newblock Geometric consistency of {M}anin's {C}onjecture.
\newblock submited, 2018.

\bibitem[LT17]{LTDuke}
B.~Lehmann and S.~Tanimoto.
\newblock On the geometry of thin exceptional sets in {M}anin's conjecture.
\newblock {\em Duke Math. J.}, 166(15):2815--2869, 2017.

\bibitem[NS20]{NS20}
M.~Nakahara and S.~Streeter.
\newblock Weak approximation and the {H}ilbert property for campana points.
\newblock submitted, 2020.

\bibitem[Pey95]{Peyre}
E.~Peyre.
\newblock Hauteurs et mesures de {T}amagawa sur les vari\'et\'es de {F}ano.
\newblock {\em Duke Math. J.}, 79(1):101--218, 1995.

\bibitem[Pey03]{Peyre03}
E.~Peyre.
\newblock Points de hauteur born\'ee, topologie ad\'elique et mesures de
  {T}amagawa.
\newblock {\em J. Th\'eor. Nombres Bordeaux}, 15(1):319--349, 2003.

\bibitem[Pey17]{Pey17}
E.~Peyre.
\newblock Libert\'{e} et accumulation.
\newblock {\em Doc. Math.}, 22:1615--1659, 2017.

\bibitem[PS20]{PS20}
M.~Pieropan and D.~Schindler.
\newblock Hyperbola method on toric varieties.
\newblock submitted, 2020.

\bibitem[PSTVA20]{PSTVA20}
M.~Pieropan, A.~Smeets, S.~Tanimoto, and A.~V\'arilly-Alvarado.
\newblock Campana points of bounded height on vector group compactifications.
\newblock {\em Proc. Lond. Math. Soc.}, 2020.
\newblock online publication.

\bibitem[Sen21]{Sen17}
A.~K. Sengupta.
\newblock Manin's conjecture and the {F}ujita invariant of finite covers.
\newblock {\em Algebra Number Theory}, 2021.
\newblock to appear.

\bibitem[ST16]{ST16}
Joseph Shalika and Yuri Tschinkel.
\newblock Height zeta functions of equivariant compactifications of unipotent
  groups.
\newblock {\em Comm. Pure Appl. Math.}, 69(4):693--733, 2016.

\bibitem[STBT07]{STBT07}
Joseph Shalika, Ramin Takloo-Bighash, and Yuri Tschinkel.
\newblock Rational points on compactifications of semi-simple groups.
\newblock {\em J. Amer. Math. Soc.}, 20(4):1135--1186, 2007.

\bibitem[Str20]{Str20}
S.~Streeter.
\newblock Campana points and powerful values of norm forms.
\newblock submitted, 2020.

\bibitem[Tat67]{Tate}
J.~T. Tate.
\newblock Fourier analysis in number fields, and {H}ecke's zeta-functions.
\newblock In {\em Algebraic {N}umber {T}heory ({P}roc. {I}nstructional {C}onf.,
  {B}righton, 1965)}, pages 305--347. Thompson, Washington, D.C., 1967.

\bibitem[TBT13]{TBT13}
Ramin Takloo-Bighash and Yuri Tschinkel.
\newblock Integral points of bounded height on compactifications of semi-simple
  groups.
\newblock {\em Amer. J. Math.}, 135(5):1433--1448, 2013.

\bibitem[VV11]{VVal11}
Karl Van~Valckenborgh.
\newblock Squareful points of bounded height.
\newblock {\em C. R. Math. Acad. Sci. Paris}, 349(11-12):603--606, 2011.

\bibitem[VV12]{VVal12}
Karl Van~Valckenborgh.
\newblock Squareful numbers in hyperplanes.
\newblock {\em Algebra Number Theory}, 6(5):1019--1041, 2012.

\bibitem[Xia20]{Xiao20}
Huan Xiao.
\newblock Campana points on biequivariant compactifications of the {H}eisenberg
  group.
\newblock submitted, 2020.

\end{thebibliography}

\end{document}